\documentclass[12pt]{article}
\long\def\kmcomment#1{}
\usepackage[leqno]{amsmath}
\usepackage{amsmath,etoolbox,mathtools}

\makeatletter
\newcommand{\leqnomode}{\tagsleft@true\let\veqno\@@leqno}
\newcommand{\reqnomode}{\tagsleft@false\let\veqno\@@eqno} 
\makeatother

\usepackage{setspace}
\usepackage{graphicx,psfrag}
\usepackage[svgnames]{xcolor}
\usepackage{tikz}
\usetikzlibrary{patterns,shapes}

\usepackage{newtxtext,newtxmath}
\usepackage[a4paper,truedimen,scale={.80,.85},centering]{geometry}

{\list{}{\leftmargin=#1 \topsep=0pt \parsep=0pt \itemsep=0pt}\item[]}%
{\endlist}
\RequirePackage{verbatim}
{\noindent\quote\endgraf\verbatim}%
{\endverbatim\endquote\endgraf\medskip}%
\makeatletter
\newenvironment{qverb*}%
{\noindent\quote\endgraf\@nameuse{verbatim*}}%
{\@nameuse{endverbatim*}\endquote\endgraf\medskip}%
\makeatother

\RequirePackage{fancyvrb} 
\DefineVerbatimEnvironment%
    {newverb}{Verbatim}
    {xleftmargin=10mm}

\DefineVerbatimEnvironment%
    {newverb*}{Verbatim}
    {xleftmargin=10mm,showspaces=true,showtabs=true} 

\usepackage{listings,jlisting,quoting}
 \RequirePackage{listings} 
 \lstnewenvironment{kmverb}
    { \lstset{xleftmargin=10mm, breaklines=true}}
    { }

\newcommand{\kmqed}{\hfill\ensuremath{\blacksquare}}

\newcommand{\futojiUp}[2]{{\mathbf #1}^{#2}}
\newcommand{\ds}{\ensuremath{\displaystyle }}
\newcommand{\pdel}{\partial}

\newcommand{\myHomW}[2]{\text{H}_{#1,#2}}

\newcommand{\mR}{\ensuremath{\mathbb{R}}} 
\newcommand{\mZ}{\ensuremath{\mathbb{Z}}} 

\newcommand{\frakg}{\mathfrak{g}}

\newcommand{\tbdl}[1]{\mathrm{T}(#1)}

\newcommand{\Sbt}[2]{[#1,#2]}                
\newcommand{\SbtS}[2]{[#1,#2]_{\text{\footnotesize S}}} 
                
\usepackage{ntheorem} 
{
\theoremheaderfont{\bfseries}
\theorembodyfont{\rmfamily}

\newtheorem{defn}{\bf Definition}
\newtheorem{prop}{Proposition}[section]
\newtheorem{exam}{Example}[section]
\newtheorem{remark}{Remark}[section]

\newtheorem{theorem}{Theorem}[section] 
\newtheorem{thm}{Theorem}[section]
\newtheorem*{thm-none}{Theorem}[section]

\newtheorem{lemma}[theorem]{Lemma}
\newtheorem{kmCor}[theorem]{Corollary}

}

\renewcommand{\[}{$$} \renewcommand{\]}{$$}

\usepackage{newtxtext} 
\usepackage[enableskew,vcentermath]{youngtab}
\usepackage{listings,jlisting,quoting}
\DefineShortVerb{\!} 

\newcommand{\myCS}[1]{ \text{C}_{#1}} 
\newcommand{\myCSW}[2]{ \text{C}_{#1,#2}} 
 
\newcommand{\mywedge}{\Delta} 
\AtBeginDocument{ \parindent=0pt }

\newcommand{\mathfrakX}[2]{\mathfrak{X}^{#1}_{#2}}
 
\newcommand{\bibunSu}[1]{\frac{\pdel }{\pdel x_{#1}}}

\newcommand{\RyakuSho}[2]{\text{Term}^{#1}_{#2}}
\newcommand{\Kai}[2]{ \xi^{#1}_{#2}}
\newcommand{\EEta}[2]{ \eta ^{#1}_{#2}}
\newcommand{\myBase}[2]{\futojiUp{\; x}{#2}\pdel_{#1}}
\newcommand{\myHTPmap}{ \pdel \circ \Phi + \phi \circ \pdel}
\newcommand{\easyXX}[2]{\mathfrak{Y}^{#1}_{#2}}

\newcommand{\kmTKu}[1]{\operatorname{UT}^{[#1]}} 
\newcommand{\kmTKd}[1]{\operatorname{DT}_{[#1]}}

\newcommand{\mySkyL}[1]{\operatorname{SkyH}_{[#1]}} 
\newcommand{\mySkyB}[1]{\operatorname{SkyB}_{[#1]}} 

\makeatletter
\newcommand\@CheckWhetherNull[1]{%
  \romannumeral0\expandafter\@secondoftwo\string{\expandafter
  \@secondoftwo\expandafter{\expandafter{\string#1}\expandafter
  \@secondoftwo\string}\expandafter\@firstoftwo\expandafter{\expandafter
  \@secondoftwo\string}\expandafter\expandafter\@firstoftwo{ }{}%
  \@secondoftwo}{\expandafter\expandafter\@firstoftwo{ }{}\@firstoftwo}%
}%
\newcommand\@CheckWhetherBlank[1]{%
  \romannumeral\expandafter\expandafter\expandafter\@secondoftwo
  \expandafter\@CheckWhetherNull\expandafter{\@firstoftwo#1{}.}%
}%
\newcommand{\kmTK}[2]{%
\operatorname{TK}^{\@CheckWhetherBlank{#1}{}{[#1]}}_{\@CheckWhetherBlank{#2}{}{[#2]}}
}%
\makeatother 
\numberwithin{equation}{section} 
\title{
The second Betti number of doubly weighted homology groups of some
pre Lie superalgebra
}

\author{
Kentaro Mikami\thanks{ 
  Department of Computer Science and Engineering 
  Akita University, partially supported by
JSPS KAKENHI Grant Number  JP26400063, JP23540067 and JP20540059.}
 \and Tadayoshi Mizutani\thanks{
Professor Emeritus, Saitama University }
}
\date{Feb 25, 2019}

\allowdisplaybreaks[1]
\begin{document}
\maketitle

\section{Introduction} 
There is a notion of (weighted) (co)homology group theory of pre Lie
superalgebras like those of Lie algebras.  
In \cite{Mik:Miz:super2}, we introduced the notion of doubly weighted
homology groups for doubly weighted (say (w,h)) pre Lie superalgebras.  
The pre Lie superalgebra we handle in this paper is 
 the exterior algebra of polynomial coefficient
multi-vector fields on n-plane with the super bracket is the Schouten
bracket. Each generator there is written as  
\[
x_{1}^{b_{1}} \cdots x_{n}^{b_{n}}
\bibunSu{a_{1}} \wedge \cdots  \bibunSu{a_{m}} \;  ( = \futojiUp{x}{\beta}
\pdel_{\alpha})\] where 
\[ \beta = (b_{1},\ldots,b_{n})\;,\;
( b_{i} \geqq 0\;\text{for}\; \forall i) \quad\text{and}\quad  
\alpha = (a_{1},\ldots,a_{m})\;,\; (1\leq a_{1} < \cdots < a_{m} \leq
n)\;. \]
Then the first weight is \(m-1\) and  
the double weight is \(( m-1, \sum_{i=1}^{n}b_{i} -1)\) by definition. 
We often denote  \(|\alpha| = m = \text{length of } \alpha\) and 
\( |\beta | =   \sum_{i=1}^{n}b_{i}\).

\medskip

In \cite{Mik:Miz:super}, we have proven a result:

(1) The Euler number is 0 for all doubly weighted homology groups. 

In \cite{Mik:Miz:super2}, we have proven several results:

(2) Each Betti number is 0 for \((w,h)\)-doubly weighted homology groups
if \(w\ne h\). 

(3) The first Betti number is 0 for all doubly weighted homology groups. 

(4) The second Betti number is 0 for doubly weighted homology groups if
\(w=h=0\). 

\medskip

In this paper, we have a main result below: 

(5) The second Betti number is 0 for  
every doubly weighted homology groups with \( w=h \).

\medskip

Thus, combining (2) in \cite{Mik:Miz:super} and (5), we conclude that 

(6) The second Betti number is 0 for all doubly weighted homology
groups. 

\section{Preliminaries}
We recall the notions and notations which we need here quickly. If it is not
enough, refer to \cite{Mik:Miz:super}. 

First we recall the definition of pre Lie
superalgebra.  
\begin{defn}[pre Lie superalgebra]
Suppose 
$\frakg $ is graded by \(\ds \mZ\) as 
\(\ds \frakg = \sum_{j\in \mZ} \frakg_{j} \)
and has a bilinear operation \(\Sbt{.}{.}\)  satisfying 
\begin{align}
& \Sbt{ \frakg_{i}}{ \frakg_{j}} \subset \frakg_{i+j} \label{cond:1} \\
& \Sbt{X}{Y} = (-1) ^{1+ x y} 
 \Sbt{Y}{X} \quad \text{ where }  X\in \frakg_{x} \text{ and } 
 Y\in \frakg_{y}  \\
& 
(-1)^{x z} \Sbt{ \Sbt{X}{Y}}{Z}  
+(-1)^{y x} \Sbt{ \Sbt{Y}{Z}}{X}  
+(-1)^{z y} \Sbt{ \Sbt{Z}{X}}{Y}  = 0 \quad \text{(Jacobi
identity).} \label{super:Jacobi}
\end{align}
Then we call \(\frakg\) a pre (or \(\mZ\)-graded) Lie superalgebra.   
\end{defn} 
In the usual Lie algebra homology theory, $m$-th chain space is the
exterior algebra \(\ds \Lambda^{m}\frakg\) of \(\frakg\) and the
boundary operator is essentially \(\ds X \wedge Y \mapsto \Sbt{X}{Y}\). 

In pre Lie superalgebras, by ``super'' skew-symmetry of bracket
operation,  $m$-th chain space \(\ds \myCS{m}\) is defined as follows:
\(\ds \myCS{m}\) is the quotient of the tensor space \(\ds \otimes^{m}
\frakg \) of \(\frakg \) by the 2-sided ideal generated by
\begin{align} & X \otimes Y + (-1)^{x y} Y \otimes X \quad\text{where
}\quad X \in \frakg_{x}, Y \in \frakg_{y} \;, \end{align} and we denote
the equivalence class of \(\ds X \otimes Y \) by \(\ds X \mywedge Y\).

Since \( \ds X_{\text{odd}} \mywedge Y_{\text{odd}}  = Y_{\text{odd}}
\mywedge X_{\text{odd}} \) and \( \ds X_{\text{even}} \mywedge
Y_{\text{any}}  = - Y_{\text{any}} \mywedge X_{\text{even}} \) hold, \(
\mywedge ^{m} \frakg_{i} \) is a symmetric algebra for odd $i$ and is a
skew-symmetric algebra for even $i$ in the usual sense. 

We introduce a recursive formula of the boundary
operator using the left action.  
\begin{align} 
\label{pdel:recurs:left}
\pdel (A_{0} \mywedge A_{1} \mywedge 
\cdots \mywedge A_{m}) &= 
- A_{0} \mywedge \pdel (A_{1} \mywedge \cdots \mywedge A_{m})
+ A_{0} \cdot  (A_{1} \mywedge \cdots \mywedge A_{m})
\\\noalign{where} \label{left:action}
A_{0} \cdot  (A_{1} \mywedge \cdots \mywedge A_{m}) &=
\Sbt{A_{0}}{A_{1}} \mywedge  
(A_{2} \mywedge \cdots \mywedge A_{m})
+ (-1)^{a_{0}a_{1}} A_{1} \mywedge \left(
A_{0} \cdot  (A_{2} \mywedge \cdots \mywedge A_{m})
\right)  \\
& = \sum_{i=1}^{m} (-1)^{ a_{0}\sum_{s<i}a_{s}} 
A_{1} \mywedge \cdots \mywedge  \Sbt{A_{0}}{A_{i}} \mywedge \cdots   
\mywedge A_{m}) 
\end{align}
for each homogeneous elements \( A_{i}\in\frakg_{a_{i}}\).  
In lower degree, the boundary operator is given
as bellows:  
\begin{align} 
\pdel ( A \mywedge B ) &= \Sbt{A}{B} \\
\pdel ( A \mywedge B \mywedge C ) &= -  {A} \mywedge {
 \Sbt{B}{C} } + 
 \Sbt{A}{B} \mywedge C + (-1)^{a b} B \mywedge  
 \Sbt{A}{C}  
\end{align}
for each homogeneous elements 
\( A\in\frakg_{a}\), \( B\in\frakg_{b}\), \( C\in\frakg_{c}\).

\begin{exam}
A prototype of pre Lie superalgebra is the exterior algebra of 
the sections of exterior power of  tangent bundle of a differentiable manifold $M$ of
dimension $n$, 
\begin{equation} \frakg = \sum_{i=1}^{n} \Lambda^{i} \tbdl{M} = 
\sum_{i=0}^{n-1} \frakg_{i} 
\;, \quad \text{where}\quad 
\frakg_{i} = \Lambda^{i+1} \tbdl{M} \;  \label{our:tbdl}
\end{equation} with the Schouten bracket. 

There are several ways defining 
the Schouten bracket, namely, axiomatic explanation, sophisticated one
using Clifford algebra or more direct ones (cf.\ 
\cite{Mik:Miz:homogPoisson}).  In the context of Lie algebra homology
theory, we introduce the Schouten bracket as follows:

\begin{defn}[Schouten bracket]
For \(\ds A \in \Lambda^{a} \tbdl{M} \) and 
\(\ds B \in \Lambda^{b} \tbdl{M} \), define a binary operation by 
\begin{equation}\label{defn:abst:schouten}
(-1)^{a+1} \SbtS{A}{B} = 
 \pdel_{0} (A\wedge B) - (\pdel_{0} A)\wedge B -
(-1)^{a} A \wedge \pdel_{0} B \;, 
\end{equation}
where \(\pdel_{0}\) is the boundary operator
in the context of Lie algebra homology of vector fields. 
\end{defn}

In some sense, the Schouten bracket measures gap of  
the boundary operator \(\pdel_{0}\)  from the derivation. 

Hereafter, we denote \( \SbtS{A}{B}\) by  \( \Sbt{A}{B}\) simply.   

\end{exam}
The first chain space is
\(\myCS{1} = \frakg= \sum_{p=1}^{n} \Lambda^{p} \tbdl{M}\).  
The second chain space is 
\begin{align*} \myCS{2} = \frakg \mywedge  \frakg  
= \sum_{1 \leq p \leq q \leq n } \Lambda^{p} \tbdl{M} \mywedge \Lambda^{q} \tbdl{M}  
=&  \Lambda^{1} \tbdl{M} \mywedge \Lambda^{1} \tbdl{M} + 
  \Lambda^{1} \tbdl{M} \mywedge \Lambda^{2} \tbdl{M} + \cdots 
\\& 
 + \Lambda^{2} \tbdl{M} \mywedge \Lambda^{2} \tbdl{M}
 + 
  \Lambda^{2} \tbdl{M} \mywedge \Lambda^{3} \tbdl{M} + \cdots 
\end{align*}

\begin{remark} \label{remark:Poisson}
Let \(\ds \pi\in\Lambda^{2}\tbdl{M}\). Then \(\ds \pi\mywedge\pi \in 
 \Lambda^{2} \tbdl{M} \mywedge \Lambda^{2} \tbdl{M} \subset \myCS{2}\)
 and \(\ds \pdel ( \pi\mywedge\pi ) = \Sbt{\pi}{\pi} \in \myCS{1}\).
 Thus, 
\(\ds \pi\in\Lambda^{2}\tbdl{M}\) is Poisson if and only if 
 \(\ds \pdel ( \pi\mywedge\pi ) = 0 \), and we express it by \(\ds \pi
 \in \sqrt{ \ker(\pdel) }\) symbolically. 
It will be interesting to study   
\(\ds \sqrt{ \ker(\pdel) }\) and also interesting to study specific
properties of Poisson structures in 
\(\ds \sqrt{ \pdel( \myCS{3} ) } \), which come from the boundary image of 
the third chain space  \(\myCS{3}\).  
\end{remark}


\subsection{First weight} 
\begin{defn}
We say a non-zero element in 
\(\ds
\frakg_{i_{1}} \mywedge  \cdots \mywedge
\frakg_{i_{m}} \) has the (first) weight  \(\ds i_{1}+\dots + i_{m} \). 
Define the subspace of 
\(\myCS{m} \) by 
\(\ds \myCSW{m}{w} = 
 \sum_{ \substack{ i_{1}\leq \ldots \leq i_{m}\\ \sum_{s=1}^{m}i_{s}  = w } } 
\frakg_{i_{1}} \mywedge  \cdots \mywedge
\frakg_{i_{m}} \), which is the direct sum of different
types of spaces but the same weight $w$.  
\end{defn}

\begin{prop}
The (first) weight $w$ is preserved by \(\ds \pdel_{ }\), i.e., we have 
\(\ds \pdel_{ }( \myCSW{m}{w} ) \subset \myCSW{m-1}{w} \). Thus, we have
for a fixed $w$, 
$w$-weighted homology groups
\[\ds 
\myHomW{m}{w} (\frakg, \mR) =  \ker(\pdel_{} : \myCSW{m}{w} \rightarrow
\myCSW{m-1}{w})/ \pdel ( \myCSW{m+1}{w} )\;.  
\] 
\end{prop}
\kmcomment{ 
In case of \eqref{our:tbdl}, possible weights are non-negative integers.  
When weight is 0,  the $m$-chain space 
is simply given by 
\(\ds \myCSW{m}{0} = \mywedge ^{m} \frakg_{0}= 
\mywedge ^{m} \tbdl{M} 
\)
and the homology is the Lie algebra homology of vector fields 
for \(m=1,\dots,n \).
For lower weights 1 or 2, 
 the chain spaces are simply given by 
\begin{align*}
\myCSW{m}{1} & = \mywedge^{m-1} \frakg_{0} \mywedge  \frakg_{1} = 
\mywedge^{m-1} \tbdl{M} \mywedge  
\Lambda^{2} \tbdl{M}\quad \text{for}\quad m=1,\dots  \;,
\\
\myCSW{m}{2} & = \mywedge^{m-1} \frakg_{0} \mywedge  \frakg_{2} 
\oplus \mywedge^{m-2} \frakg_{0} \mywedge^{2}\frakg_{1} 
\\& 
= 
\mywedge^{m-1} \tbdl{M} \mywedge  \Lambda^{3} \tbdl{M}
\oplus \mywedge^{m-2} \tbdl{M} \mywedge^{2}\Lambda^{2} \tbdl{M} 
\quad \text{for}\quad m=1,\dots  \;.
\end{align*} 
\begin{remark}
In particular, 
\(\ds \myCSW{1}{2}  = \Lambda^{3} \tbdl{M}\), 
\(\ds \myCSW{2}{2}  =
\tbdl{M} \mywedge  \Lambda^{3} \tbdl{M}
\oplus \Lambda^{2} \tbdl{M} \mywedge\Lambda^{2} \tbdl{M} \),  
\(\ds \myCSW{3}{2}  =
\tbdl{M} \mywedge  \tbdl{M} \mywedge  
\Lambda^{3} \tbdl{M}
\oplus \tbdl{M} \mywedge\Lambda^{2} \tbdl{M}  \mywedge\Lambda^{2} \tbdl{M} 
\).  Thus, by introducing weight, the chain spaces become smaller and
research becomes a little clear and easier.   
\end{remark}
}

\subsection{Double weight} 
\begin{defn}[Double-weight] \label{defn:w:weight}
Assume that each subspace \(\ds \frakg_{i} \) of a given pre Lie superalgebra
\(\frakg\) is direct decomposed by subspaces \( \ds \frakg_{i,j}\)
as \(\ds \frakg_{i}  = \sum_{j} \frakg_{i,j} \) and satisfies 
\begin{equation}
\Sbt{ X }{ Y } \in \frakg_{i_1+i_2, j_1+j_2} \quad \text{for each}\ 
 X\in \frakg_{i_1,j_1} \ ,  \ 
Y\in \frakg_{i_2,j_2} \;. \end{equation}  
We say such pre Lie superalgebras are double-weighted. 

We may define double-weighted $m$-th chain space by
\[
\myCSW{m}{w,h} = \sum_{
\substack{ 
i_{1}\leq \ldots \leq i_{m}\;,\; 
\sum_{s=1}^{m}i_{s}  = w \\
\sum_{s=1}^{m}h_{s}  = h 
} } 
\frakg_{i_{1},h_{1}} \mywedge  \cdots \mywedge
\frakg_{i_{m},h_{m}} 
\]
\end{defn}

\begin{prop}
The double-weight $(w,h)$ is preserved by \(\ds \pdel_{ }\), i.e., we have 
\(\ds \pdel_{ }( \myCSW{m}{w,h} ) \subset \myCSW{m-1}{w,h} \). Thus, we have 
$(w,h)$-weighted homology groups
\[\ds 
\myHomW{m}{w,h} (\frakg, \mR) =  \ker(\pdel_{} : \myCSW{m}{w,h} \rightarrow
\myCSW{m-1}{w,h})/ \pdel ( \myCSW{m+1}{w,h} )\;.  
\] 
\end{prop}

As we explained in Introduction, 
we consider the Euclidean space \(\ds M = \mR^
{n} \) with the Cartesian coordinates \(x_{1},\dots,x_{n}\). 
Then, we get a pre Lie super subalgebra consisting of multi vector fields of
polynomial coefficients. 
We define \begin{align*} \frakg_{i,j} & =
\mathfrakX{i+1}{j+1} (\mR^{n}) 
 = \{
(i+1)\text{-multi vector fields with }
(j+1)\text{-homogeneous polynomials}\} \;. \end{align*} 
We see easily that \(\ds \Sbt{ \frakg_{i_1,j_1}} 
{ \frakg_{i_2,j_2}} \subset \frakg_{i_1+i_2, j_1+j_2} \) and so we get a
double-weighted pre Lie superalgebra.

\section{2nd Betti number for  general \(w=h\)} 
2-chain space \(\myCSW{2}{w,h}\) consists of 
\(\ds \mathfrakX{a_{1}}{b_{1}}\mywedge \mathfrakX{a_{2}}{b_{2}}\) with 
            \( a_{1} +a_{2} = 2+w\) and \( b_{1} +b_{2} = 2+w\). 
Since \(\mywedge\) has a super symmetric property, 
\(\ds \mathfrakX{a_{1}}{b_{1}}\mywedge \mathfrakX{a_{2}}{b_{2}} 
= \mathfrakX{a_{2}}{b_{2}}\mywedge \mathfrakX{a_{1}}{b_{1}}
\)  as spaces,
so we may express each subspace 
\(\ds \mathfrakX{a_{1}}{b_{1}}\mywedge \mathfrakX{a_{2}}{b_{2}} 
= \mathfrakX{a_{2}}{b_{2}}\mywedge \mathfrakX{a_{1}}{b_{1}} \)  
by \(\ds \easyXX{a_{1}}{b_{1}} := 
\mathfrakX{a_{1}} {b_{1}} \mywedge 
\mathfrakX{2+w-a_{1}} {2+w-b_{1}} \) with \( 1\leqq a_{1} \leqq 1+ w/2\).

\(\ds \myCSW{2}{w,w}\) has two expressions depending on \(w\)'s
parity.    
\begin{align*}\text{If \(w =2\Omega +1\), then\quad} 
\myCSW{2}{w,w} &=  \sum_{ a_{1}=1}^{\Omega +1}\sum_{ b_{1}=0}^{2+w} 
\mathfrakX{a_{1}}{b_{1}} \mywedge \mathfrakX{2+w - a_{1}}{2+w -b_{1}} \\
\text{If \(w =2\Omega \), then\quad} 
\myCSW{2}{w,w} &=  \sum_{ a_{1}=1}^{\Omega }\sum_{ b_{1}=0}^{2+w} 
\mathfrakX{a_{1}}{b_{1}} \mywedge \mathfrakX{2+w - a_{1}}{2+w -b_{1}} 
+ \sum_{ b_{1}=0}^{\Omega } \mathfrakX{\Omega + 1}{b_{1}} 
 \mywedge \mathfrakX{\Omega + 1}{2+w -b_{1}} 
+ \mathfrakX{\Omega +1}{\Omega +1} \mywedge \mathfrakX{\Omega +1}{\Omega+1} 
\end{align*}

To understand how \(\ds\myCSW{2}{w,w}\) is decomposed, we plot the point 
 \( (a_{1},b_{1})\) of \(\easyXX{a_{1}}{b_{1}}\)    
on the 2-plane.  The next are examples of \(w=3,4\).

\hspace*{30mm} 
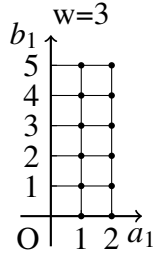
\begin{minipage}[t]{0.38\textwidth}
\begin{tikzpicture}[scale = 0.4]
\draw[thick, ->] (-1,0) -- (3,0) node [below] {$a_{1}$};
\draw[thick, ->] (0,-1) -- (0,6) node [left] {$b_{1}$};
\node at (0,0) [anchor=north east] {O};
\node at (1,6) [anchor=south ] {w=3};
\foreach \x in {1,2}
   \draw (\x cm,.2) -- (\x cm,0) node[anchor=north] {$\x$};
\foreach \y in {1,2,3,4,5}
   \draw (.2,\y cm) -- (0,\y cm) node[anchor=east] {$\y$};
\draw (0,0) grid  (2,5); 
\foreach \x in {4,5} \fill (2,\x) circle (3pt) ;
\fill (1,5) circle (3pt) ;
\foreach \x in {0,1,2,3,4} \fill (1,\x) circle (3pt) ;
\foreach \x in {0,1,2,3} \fill (2,\x) circle (3pt) ;
\end{tikzpicture}
\end{minipage} 
\begin{minipage}[t]{0.52\textwidth}
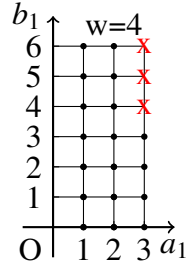
\begin{tikzpicture}[scale = 0.4]
\draw (0,0) grid  (3,6); 
\draw[thick, ->] (-1,0) -- (4,0) node [below] {$a_{1}$};
\draw[thick, ->] (0,-1) -- (0,7) node [left] {$b_{1}$};
\node at (0,0) [anchor=north east] {O};
\node at (2,6) [anchor=south ] {w=4};
\foreach \x in {1,2,3}
   \draw (\x cm,.2) -- (\x cm,0) node[anchor=north] {$\x$};
\foreach \y in {1,2,3,4,5,6}
   \draw (.2,\y cm) -- (0,\y cm) node[anchor=east] {$\y$};
\foreach \x in {1,2}
\foreach \y in {0,1,2} \fill (\x,\y) circle (3pt) ;

\foreach \y in {4,5,6} \fill[color=red] (3,\y) node {x} ;

\foreach \y in {4,5} \fill  (1,\y) circle (3pt) ; 
\foreach \y in {6} \fill (1,\y) circle (3pt) ; 

\foreach \y in {4} \fill  (2,\y) circle (3pt) ; 
\foreach \y in {5,6} \fill  (2,\y) circle (3pt) ; 

\foreach \x in {1,2,3} \foreach \y in {0,1,2,3} 
\fill (\x,\y) circle (3pt) ;

\end{tikzpicture}
\end{minipage}

\begin{defn}
We define the type of 
single 2-chain 
\(\ds \myBase{A_{1}}{B_{1}} \mywedge \myBase{A_{2}}{B_{2}} \) by 
\begin{align*} 
\text{type \textit{TR}} & \quad  
\text{if}\quad 
 |A_{1}|+|B_{1}| < |A_{2}|+|B_{2}| \text{\quad or\quad}  
 |A_{1}|+|B_{1}| = |A_{2}|+|B_{2}| \text{\quad \&\quad}  
 |A_{1}| \leqq  |A_{2}| \\
\text{ type  \textit{TL}} &\quad   \text{otherwise}.  
 \end{align*} 
For each 
\(\ds \myBase{A_{1}}{B_{1}} \mywedge \myBase{A_{2}}{B_{2}} \), 
we define a 3-chain by     
\begin{align}
&  \sum_{\ell=1}^{n} 
\pdel_{\ell}\mywedge(\myBase{A_{1}}{B_{1}})\mywedge
(x_{\ell}\myBase{A_{2}}{B_{2}}) \quad \text{ if type \it{TR} } \\
& \sum_{\ell=1}^{n} 
\pdel_{\ell} \mywedge (x_{\ell}
\myBase{A_{1}}{B_{1}})\mywedge (\myBase{A_{2}}{B_{2}})
 \quad \text{ if type \it{TL} } 
\end{align}
and extend linearly on \(\ds\myCSW{2}{w,w}\) and  denote it by 
\(\Phi\).
Define a linear map \(\ds\phi : \myCSW{1}{w,w} \to  \myCSW{2}{w,w} \) by  
\begin{equation}\phi(U)
= \sum_{\ell=1}^{n} \pdel _{\ell} \mywedge ( x_{\ell} U )\;. \end{equation}

\end{defn}

\begin{lemma}
Consider 
\(\ds \myBase{A_{1}}{B_{1}} \mywedge \myBase{A_{2}}{B_{2}} \). 
\begin{align*} 
\shortintertext{If \textit{TR}, then}   
(\pdel \Phi + \phi \pdel)  
(\myBase{A_{1}}{B_{1}} \mywedge \myBase{A_{2}}{B_{2}}) 
&=  - \sum_{\ell=1}^{n} \pdel_{\ell} \mywedge 
\Sbt{\futojiUp{x}{B_{1}} \pdel_{A_{1}} }{ x_{\ell}}
\futojiUp{x}{B_{2}} \pdel_{A_{2}}   
 + \sum_{\ell=1}^{n} \Sbt{\pdel_{\ell} }
{\futojiUp{x}{B_{1}} \pdel_{A_{1}} }\mywedge( x_{\ell}
\futojiUp{x}{B_{2}} \pdel_{A_{2}}) 
\\&\quad 
 + (n + |B_{2}|)
(\futojiUp{x}{B_{1}} \pdel_{A_{1}} )\mywedge
( \futojiUp{x}{B_{2}} \pdel_{A_{2}} )  \;.
\shortintertext{If \textit{TL}, then}   
(\pdel \Phi + \phi \pdel)  
(\myBase{A_{1}}{B_{1}} \mywedge \myBase{A_{2}}{B_{2}}) 
&=  
- (-1)^{|A_{1}|(|A_{2}|+1)} \sum_{\ell=1}^{n} \pdel_{\ell} \mywedge  
\Sbt { x_{\ell}} {\futojiUp{x}{B_{2}}\pdel_{A_{2}}} 
\futojiUp{x}{B_{1}} \pdel_{A_{1}} 
\\ 
&\quad + (n +| B_{1}|)
(\futojiUp{x}{B_{1}} \pdel_{A_{1}} )\mywedge
( \futojiUp{x}{B_{2}} \pdel_{A_{2}} )  
 + \sum_{\ell=1}^{n} 
( x_{\ell}\futojiUp{x}{B_{1}} \pdel_{A_{1}} )\mywedge
\Sbt{ \pdel_{\ell} } { \futojiUp{x}{B_{2}} \pdel_{A_{2}}}\;.   
\end{align*}

For an element in \( \easyXX{a}{b} = \mathfrakX {a}{b} \mywedge
\mathfrakX {2+w-a}{2+w-b}\), when we want only to indicate its type, we
sometimes use the symbol \(\RyakuSho{a}{b}\) for this purpose.  Using
this abbreviation and denoting 
\(\pdel \Phi + \phi \pdel\)  by \(\Psi\), the above become more simpler. 
\begin{align} 
\shortintertext{If \textit{TR}, then}   
\Psi   
(\myBase{A_{1}}{B_{1}} \mywedge \myBase{A_{2}}{B_{2}}) 
&=  
  (n + |B_{2}|) (\futojiUp{x}{B_{1}} \pdel_{A_{1}} )\mywedge
( \futojiUp{x}{B_{2}} \pdel_{A_{2}} ) + 
\RyakuSho{|A_{1}|}{|B_{1}|-1} + \RyakuSho{1}{0}   \label{map::DN}
\shortintertext{If \textit{TL}, then}   
\Psi 
(\myBase{A_{1}}{B_{1}} \mywedge \myBase{A_{2}}{B_{2}}) 
&=  
 (n + |B_{1}|)
(\futojiUp{x}{B_{1}} \pdel_{A_{1}} )\mywedge
( \futojiUp{x}{B_{2}} \pdel_{A_{2}} )  
+ 
\RyakuSho{|A_{1}|}{|B_{1}|+1} + \RyakuSho{1}{0}  \label{map::UP}
\end{align}

\end{lemma}
\textbf{Proof:}
When \textit{TR}, 
\begin{align*}
&\quad 
\pdel( \Phi 
(\myBase{A_{1}}{B_{1}} \mywedge \myBase{A_{2}}{B_{2}}) ) 
= \pdel (
\sum_{\ell=1}^{n} \pdel_{\ell} \mywedge 
(\myBase{A_{1}}{B_{1}}) \mywedge (x_{\ell} \myBase{A_{2}}{B_{2}}) ) 
\\&= 
- \sum_{\ell=1}^{n} \pdel_{\ell} \mywedge  
\Sbt{\futojiUp{x}{B_{1}} \pdel_{A_{1}} }{ x_{\ell}
\futojiUp{x}{B_{2}} \pdel_{A_{2}}}  
+ \sum_{\ell=1}^{n} 
 \Sbt{\pdel_{\ell} }
{\futojiUp{x}{B_{1}} \pdel_{A_{1}} }\mywedge( x_{\ell}
\futojiUp{x}{B_{2}} \pdel_{A_{2}}) 
+ \sum_{\ell=1}^{n}
(\futojiUp{x}{B_{1}} \pdel_{A_{1}} )\mywedge
\Sbt{ \pdel_{\ell} } { x_{\ell}
\futojiUp{x}{B_{2}} \pdel_{A_{2}}}  \\
&= 
- \phi( \pdel 
(\myBase{A_{1}}{B_{1}} \mywedge \myBase{A_{2}}{B_{2}}) 
) 
- \sum_{\ell=1}^{n} \pdel_{\ell} \mywedge
\Sbt{\futojiUp{x}{B_{1}} \pdel_{A_{1}} }{ x_{\ell}}
\futojiUp{x}{B_{2}} \pdel_{A_{2}}  \\ 
&\quad + \sum_{\ell=1}^{n} 
 \Sbt{\pdel_{\ell} }
{\futojiUp{x}{B_{1}} \pdel_{A_{1}} }\mywedge( x_{\ell}
\futojiUp{x}{B_{2}} \pdel_{A_{2}}) 
 + (n+ |B_{2}|)
(\futojiUp{x}{B_{1}} \pdel_{A_{1}} )\mywedge
( \futojiUp{x}{B_{2}} \pdel_{A_{2}} )  \;.
\end{align*}

When \textit{TL}, 
\begin{align*}
&\quad 
\pdel( \Phi 
(\myBase{A_{1}}{B_{1}} \mywedge \myBase{A_{2}}{B_{2}}) ) 
= \pdel (
\sum_{\ell=1}^{n} \pdel_{\ell} \mywedge 
(x_{\ell} \myBase{A_{1}}{B_{1}}) \mywedge (\myBase{A_{2}}{B_{2}}) ) 
\\
&= - \sum_{\ell=1}^{n} \pdel_{\ell} \mywedge  
\Sbt{ x_{\ell}\futojiUp{x}{B_{1}} \pdel_{A_{1}} }{
\futojiUp{x}{B_{2}} \pdel_{A_{2}}}  
+ \sum_{\ell=1}^{n} 
 \Sbt{\pdel_{\ell} }
{ x_{\ell}\futojiUp{x}{B_{1}} \pdel_{A_{1}} }\mywedge(
\futojiUp{x}{B_{2}} \pdel_{A_{2}}) 
+ \sum_{\ell=1}^{n}
( x_{\ell}\futojiUp{x}{B_{1}} \pdel_{A_{1}} )\mywedge
\Sbt{ \pdel_{\ell} } { \futojiUp{x}{B_{2}} \pdel_{A_{2}}}  \\
&= 
- \phi( \pdel 
(\myBase{A_{1}}{B_{1}} \mywedge \myBase{A_{2}}{B_{2}}) ) ) 
 - (-1)^{|A_{1}|(|A_{2}|+1)} \sum_{\ell=1}^{n} \pdel_{\ell} \mywedge  
\Sbt { x_{\ell}} {
\futojiUp{x}{B_{2}} \pdel_{A_{2}} } 
\futojiUp{x}{B_{1}} \pdel_{A_{1}} 
\\ &\quad  
+  (n +|B_{1}|)
(\futojiUp{x}{B_{1}} \pdel_{A_{1}} )\mywedge
( \futojiUp{x}{B_{2}} \pdel_{A_{2}} ) 
+ \sum_{\ell=1}^{n} 
( x_{\ell}\futojiUp{x}{B_{1}} \pdel_{A_{1}} )\mywedge
\Sbt{ \pdel_{\ell} } { \futojiUp{x}{B_{2}} \pdel_{A_{2}}}  \;. 
\end{align*}
\kmqed

\begin{remark}
\eqref{map::DN} and  \eqref{map::UP} say that  
\( \myBase{A_{1}}{B_{1}} \mywedge \myBase{A_{2}}{B_{2}} \) 
is one step  descending  to \(\RyakuSho{|A_{1}|}{|B_{1}|-1}\) 
if \textit{TR} type and one step   
ascending to \(\RyakuSho{|A_{1}|}{|B_{1}|+1}\) if \textit{TL} type by the map
\( \Psi := \pdel\circ \Phi + \phi\circ \pdel \) modulo \(\RyakuSho{1}{0}\).  
Hereafter, sometimes or somewhere we express \(\Psi(U)\) by 
omitting \(\ds \RyakuSho{1}{0}\) but never forget its contribution. 
Actually, we see next interesting property. 
\end{remark}
\begin{lemma} \label{lemma:interest}
If \( U \in \easyXX{1}{0}\) then it satisfies \(\Psi( U ) = (n+w+1) U\).  
\end{lemma}
\textbf{Proof:}
We may express 
\(\ds 
U = \sum_{\substack{i, \beta, G^{i,\beta}}} \pdel_{i}
\mywedge G^{i,\beta} \pdel_{\beta}\) with \(|\beta|=1+w, |G^{i,\beta}|=2+w\). 
\begin{align*}
\pdel \Phi( U) &=
\sum_{\ell} \sum_{\substack{i, \beta, G^{i,\beta}}}
\pdel ( \pdel_{\ell} \mywedge \pdel_{i} \mywedge x_{\ell} G^{i,\beta} \pdel_{\beta}) 
= 
\sum_{\ell} \sum_{\substack{i, \beta , G^{i,\beta}}}
( - \pdel_{\ell} \mywedge  \Sbt{ \pdel_{i}} {x_{\ell} G^{i,\beta} \pdel_{\beta} }
 +  \pdel_{i} \mywedge \Sbt{ \pdel_{\ell}}{ x_{\ell} G^{i,\beta} \pdel_{\beta}} )
\\ & 
= - \phi \pdel U - 
\sum_{\ell} \sum_{\substack{i, \beta,  G^{i,\beta}}}
\pdel_{\ell} \mywedge  \Sbt{ \pdel_{i}} {x_{\ell}} G^{i,\beta} \pdel_{\beta} 
 +  (n+w+2)  U  
= - \phi \pdel U  +  (n+w+1)  U \;,  \\
\Psi(  U ) &=   (n+w+1)  U \;.  
\end{align*}
\kmqed

To know the type of 
\(\ds \easyXX{a_{1}}{b_{1}} \) is easy as follows:
\begin{align*}
\text{
\textit{TL} type iff }\ a_{1}+ b_{1} > 2+w \;,\; \text{we express a general element there by  }\ 
 \Kai{a_{1}}{b_{1}}\;. \\
 \text{
\textit{TR} type iff }\ a_{1}+ b_{1} \leqq  2+w\;,\; \text{we
express a general element there by  }\ 
 \EEta{a_{1}}{b_{1}}\;. \end{align*}

We divide 
\(\ds \myCSW{2}{w,w}  \) into 
\kmcomment{
=  \sum_{a_{1}+b_{1} \leqq 2+w} \easyXX{a_{1}}{b_{1}}+
\sum_{a_{1},b_{1}>2+w} \easyXX{a_{1}}{b_{1}} \)  
 consists of  } the two subspaces
\begin{alignat*}{2}
  W_{[TR]}
&:=\sum_{a_{1}+b_{1} \leqq 2+w} \easyXX{a_{1}}{b_{1}}
\;\text{  of type}\; \textit{TR}\;,
& \quad 
\text{and}\quad   
  W_{[TL]} & := 
\sum_{a_{1} + b_{1}>2+w} \easyXX{a_{1}}{b_{1}} \;\text{ 
of type \textit{TL}}\; . 
\end{alignat*}

\begin{prop}
The subspace \( W_{[TR]}\) is invariant under \(\ds \Psi (= \myHTPmap)\).  

The subspace \( W_{[TL]}\) is invariant under \(\ds \Psi (= \myHTPmap)\)
modulo   \(\RyakuSho{1}{0}\).  
\end{prop}

\begin{defn} 
Let 
\begin{alignat*}{3}
\kmTK{\ell}{} &:= 
\sum_{\mathclap{a_{1}+b_{1} = \ell + 2 +w}} \easyXX{a_{1}}{b_{1}} \;, 
\quad &
\kmTK{\ell,s}{} & := 
\sum_{\mathclap{\substack{a_{1}+b_{1} = \ell + 2 +w\\ a_{1} \geqq s}}} 
\easyXX{a_{1}}{b_{1}} \subset \kmTK{\ell}{} \;,  
\quad & 
\kmTKu{\ell} &:= 
\sum_{\mathclap{a_{1}+b_{1} > \ell + 2 +w}} \easyXX{a_{1}}{b_{1}} \;. 
\end{alignat*}  
\end{defn}
We easily see

\begin{prop} 
\(\ds \kmTKu{\ell} = \kmTK{\ell+1}{}+\kmTKu{\ell+1}\) for \(\ell\geqq
0\) and  
\(\ds W_{[TL]} = \kmTKu{0}\). 
\begin{align*}
&  W_{[TL]} 
= \sum_{\ell=1}^{\Omega_{e} } \kmTK{\ell}{}\;,\quad \text{ where  }\;
\Omega_{e} := \begin{cases} 
\Omega +1 & \ \text{ if }\ w=2\Omega +1 \\
\Omega  & \ \text{ if }\ w=2\Omega  \end{cases}\;.\\
& 
 \Psi( \kmTKu{\ell}) = \kmTKu{\ell}\;,\quad(\text{in precise}\;  
 \Psi( \kmTKu{\ell}) \subset \kmTKu{\ell} + \easyXX{1}{0} )\;.\\   
& 
\Psi( \kmTK{\ell}{}) = \kmTK{\ell}{} + \kmTK{\ell +1}{}\;,\quad(\text{in precise}\;  
\Psi( \kmTK{\ell}{}) \subset \kmTK{\ell}{} + \kmTK{\ell +1}{}
+ \easyXX{1}{0})
\;.
\end{align*}
\end{prop}

Our main result is the following:
\begin{thm} The second Betti number is zero for \(\ds\{ \myCSW{\bullet}{w,w}\}\). 
\end{thm}
In order to prove the theorem above, we follow three steps:
\begin{enumerate}
\item
We reduce our discussion from \(\ds W_{[TL]} + W_{[TR]}\)
to \(\ds W_{[TR]}\).

\item
We decompose \(\ds W_{[TR]} = \sum_{a_{1}+b_{1} \leqq 2 +w, b_{1}>
1+\Omega_{e}} \easyXX{a_{1}}{b_{1}}  + \sum_{a_{1}+b_{1} \leqq 2 +w,
b_{1}\leqq 1+\Omega_{e}} \easyXX{a_{1}}{b_{1}}\).  

We reduce our
discussion to the rectangular region \(\ds \sum_{a_{1}+b_{1} \leqq 2 +w,
b_{1}\leqq 1+\Omega_{e}} \easyXX{a_{1}}{b_{1}}\).  

\item
We finish our discussion on the rectangular region. 
\end{enumerate}

\subsection{\textit{TL}}
We prepare one of two key lemmas:
\begin{lemma} \label{lemma:key}
Let \(a\) and \(s\) satisfy with \( 1 \leqq  a \leqq s\leqq \Omega_{e}\).   
\begin{align}
\shortintertext{Take}
 U & = u + u' 
\quad \text{where}\quad  u \in\kmTK{a,s}{}\smallsetminus \kmTK{a,s+1}{}  
\;\text{and}\;    
u' \in \kmTKu{a} + W_{[TR]}   \;.  
\shortintertext{Then} 
\Psi( U ) & = c U + U' \quad \exists c \ne 0\;,\;  
 U' \in \kmTK{a, \color{red}1+s}{}+ 
 \kmTKu{a} + W_{[TR]}\;. 
\end{align} 
\end{lemma}
\textbf{Proof:} 
In general, we handle an element \(U\) belong  to 
\(\ds \kmTK{a}{}+\kmTKu{a} +  W_{TR}\), 
usually we write as follows:
\( \ds U = u_{1} + u_{2} + u_{3}\in  \kmTK{a}{}+\kmTKu{a} +  W_{TR}\)
or 
\( \ds U = u_{1} + u_{2} + u_{3}\) 
where   \( \ds u_{1} \in  \kmTK{1}{}\),  
\( \ds u_{2} \in  \kmTKu{1} \),  
\( \ds u_{3} \in  W_{[TR]}\).  
Sometimes or somewhere in this article,  we impolitely write  
\( \ds U = u_{1} + \kmTKu{a} +  W_{[TR]}\) where 
\( \ds u_{1} \in  \kmTK{a}{}\) for instance when we do not refer 
\(u_{2}\) and \(u_{3}\) later in precise.    
\begin{align*}
\shortintertext{Take}
 U & = \sum_{s \leqq t \leqq \Omega_{e} } \Kai{t}{w+2+a-t}
+ u'\quad \text{where}\quad  \Kai{s}{w+2+a-s} \ne 0 \;. 
\shortintertext{Then} 
\Psi(U) &=\Psi(
\sum_{t=s}^{\Omega_{e}} \Kai{t}{ w+2+a - t } ) 
+  \Psi(  u' )  
= \sum_{t=s}^{\Omega_{e}}( (n+w+2+a-t)  \Kai{t}{w+2+a - t } +
\RyakuSho{t}{w+3+a-t} )
+ \Psi(  u' )  \\
&= \sum_{t=s}^{\Omega_{e}} (n+w+2+a-t)  \Kai{t}{ w+2+a - t } +
   \kmTKu{a} + \Psi(  u' )  \\ 
&= (n+w+2+a-s) \sum_{t=s}^{\Omega_{e}} \Kai{t}{w+2+a  - t } +
\sum_{t=s}^{\Omega_{e}} (s-t)  \Kai{t}{ w+2+a - t }+\kmTKu{a}+\Psi(u')\\ 
&= (n+w+2+a-s) (\sum_{s=t}^{\Omega_{e}} \Kai{t}{ w+2+a - t }+ u') +
\sum_{t=s}^{\Omega_{e}} (s-t)  \Kai{t}{ w+2+a - t }
\\&\quad 
-(n+w+2+a-s) u' +\kmTKu{a}+\Psi(u')\\ 
&= (n+w+2+a-s) U  +
\sum_{t=1+s}^{\Omega_{e}}(s- t)  \Kai{t}{ w+a+2 - t }
+\kmTKu{a}+ W_{[TR]}\\ 
&= (n+w+2+a-s) U  + U'  
\shortintertext{where }
U'  &= 
\sum_{t=s+1}^{\Omega_{e}}(s- t)  \Kai{t}{ w+2+a - t }
+\kmTKu{a}+ W_{[TR]}\; . 
\end{align*}   
\kmqed
\begin{kmCor} \label{cor:ascending}
Let \(U \ne 0\) and \(U \in\kmTK{a}{}\). 
Then we have some non-zero numbers \(\ds \{ c_{i} \}\) so that 
\[ (\Psi+c_{1}) \circ \cdots \circ (\Psi+c_{\ell}) U \in \kmTKu{a}\;.\] 
\end{kmCor}

\begin{prop} \label{prop:TL:to:TR}
Take  \(\ds U \in \myCSW{2}{w,w} =  W_{[TL]} + W_{[TR]}\).  
We reduce \(U\)  to \(W_{[TR]}\) in the
following sense.  
\[ (\Psi+c_{1}) \circ \cdots \circ (\Psi+c_{\ell}) U \in  W_{[TR]} \;,\] 
for some non-zero numbers \(\ds \{ c_{i} \}\).  
\end{prop}
\textbf{Proof:} 
\( U = U_{0} + U_{1}\) with \(U_{0}\in W_{[TL]}, U_{1} \in W_{[TR]}\).  
We may assume \(U_{0} \ne 0\).  
Apply Proposition \ref{cor:ascending}
several times,  we see that 
\(\ds (\Psi+c_{1}) \circ \cdots \circ (\Psi+c_{\ell}) U_{0} \in
W_{[TR]}\)  
for some non-zero numbers \(\ds \{ c_{i} \}\), and   
\(\ds  (\Psi+c_{1}) \circ \cdots \circ (\Psi+c_{\ell}) U_{1} \in  W_{[TR]} \)
because  \(W_{[TR]}\) is invariant under the action of \(\Psi\). Thus,
we see 
\(\ds (\Psi+c_{1}) \circ \cdots \circ (\Psi+c_{\ell}) U \in  W_{[TR]} \).  
  \kmqed

\subsection{\textit{TR}} 
Here we study  \(\Psi(U)\) of \(U \in W_{[TR]}\). For that purpose, we
divide  \(\ds W_{[TR]}\) into two parts, one is roof part \(\ds
\sum_{\substack{a_{1}+b_{1} \leqq 2 +w\\ b_{1}> 1+\Omega_{e}}} \easyXX{a_{1}}{b_{1}}
\)  and the other is rectangular basic part \(\ds
\sum_{\substack{a_{1}+b_{1} \leqq 2 +w \\
b_{1}\leqq 1+\Omega_{e}}} \easyXX{a_{1}}{b_{1}}\).  
We define three subspaces: 
\begin{alignat*}{3} 
\kmTK{}{p} & := \sum_{a+b=p} \easyXX{a}{b}
\;,\quad 
& \kmTK{}{p,s} & := \sum_{a+b=p, a\leqq s} \easyXX{a}{b} 
\;,\quad 
&
\kmTKd{p} &:= \sum_{a+b< p} \easyXX{a}{b}\;. 
\end{alignat*}
\begin{prop}
\( W_{TR} = \kmTKd{2+w+1} = \sum_{ p \leqq 2+w}  \kmTK{}{p} \)\;,
\quad 
\(\kmTKd{p+1} = \kmTK{}{p} + \kmTKd{p}\)\;, 

\(\Psi(\kmTK{}{p})=\kmTK{}{p} + \kmTK{}{p-1}\)\;, \quad 
\(\Psi(\kmTKd{p})= \kmTKd{p}\)\;. 
\end{prop}

Making use of descending property of \textit{TR} for \(\Psi\), 
we prepare another key lemma  like Lemma \ref{lemma:key}. 
\begin{lemma} \label{lemma:key2}
Take \(\ds  U  \in \kmTK{}{p,s}\smallsetminus  \kmTK{}{p,s-1} \). 
Then  
\(\ds (\Psi + c) U \in \kmTK{}{p,s-1} + \kmTKd{p}\) 
for some \(c \ne 0\).  
\end{lemma}
\textbf{Proof:}
\begin{align*} U &= \sum_{a+b=p, a\leqq s} \EEta{a}{b}\;,   \quad  
(\EEta{s}{p-s} \ne 0)\;.  \\
\Psi( U ) &= \sum_{a+b=p, a\leqq s}\Psi( \EEta{a}{b}) 
= \sum_{a\leqq s}( (n+2+w-p+a)\EEta{a}{p-a}+
\RyakuSho{a}{p-a-1}) \\
& =  (n+2+w-p+s)\sum_{a\leqq s} \EEta{a}{p-a}
+ \sum_{ a\leqq s} (a - s)\EEta{a}{p-a}
+ \kmTKd{p}\;  \\
& =  (n+2+w-p+s) U + \sum_{ a\leqq s} (a - s)\EEta{a}{p-a}
+ \kmTKd{p}\;.  
\end{align*} 
\kmqed

\begin{prop} \label{prop:TR:to:Rec}
Take a \(\ds U\in  W_{[TR]}\).  
We reduce \(U\)  to the rectangular part in the following sense.  
\[
(\Psi + c_{1}) \circ \cdots \circ
(\Psi + c_{k}) U \in \text{Rectangular part of }\; W_{[TR]} \;   
\] 
for some non-zero numbers \(\ds \{ c_{i} \}\).   
\end{prop}
\textbf{Proof:} 
Making use of descending property of \textit{TR}, 
applying Lemma \ref{lemma:key2}
several times along   
\(\ds \sum_{a_{1}+b_{1} = p } \easyXX{a_{1}}{b_{1}} \)  
for \(p= w+2,\ldots, 2+\Omega_{e}\). 

\kmqed

We divide the rectangular region horizontally, i.e., 
\begin{align*}
\mySkyL{ b } &:= \sum_{1 \leqq a \leqq \Omega+1} \easyXX{a}{b} \;,\quad 
\mySkyB{ b } := \sum_{b' < b} \mySkyL{ b' } \;.
\end{align*}
\(\ds \Psi( \mySkyL{ b } ) \subset \mySkyL{ b }+  \mySkyL{ b-1 } \)  
 and  \(\ds \mySkyB{ b }\) is invariant by \( \Psi = \myHTPmap \).  
We have an easy lemma.

\begin{lemma}
Take \(\ds U\in\mySkyB{ b +1 } \smallsetminus  \mySkyB{ b } \). 
Then   
\[ \Psi(U) = (n+w+2-b) U + U' \quad\text{where}\ U' \in 
\mySkyB{ b }\;. \]
\end{lemma}
\textbf{Proof:}
\(\ds U = \sum_{a=1}^{\Omega+1} \EEta{a}{b}\) and   
\begin{align*}
\Psi(U) &= \sum_{a=1}^{\Omega+1} \Psi(\EEta{a}{b} ) 
 = \sum_{a=1}^{\Omega+1}((n+w+2-b)\EEta{a}{b} + \RyakuSho{a}{b-1}) 
 = (n+w+2-b)\sum_{a=1}^{\Omega+1}\EEta{a}{b} + \mySkyB{b}\\
 &= (n+w+2-b) U + U'\;, \quad U' \in \mySkyB{b}\;. 
 \end{align*}
\kmqed

\begin{prop} \label{prop:Rec:to:Fin}
Take \(U\) from the rectangle region of \(W_{[TR]}\) defined by
\(\ds \sum_{\substack{ 1\leq a_{1} \leq \Omega +1 \\ 0 \leqq b_{1} \leqq  
1+\Omega_{e}}}\easyXX{a_{1}}{b_{1}}\).  Then 
we reduce \(U\)  to the null in the 
following sense.  
\[ (\Psi + c_{1}) \circ \cdots \circ (\Psi + c_{k+1}) U  = 0 \] 
for some non-zero numbers \(\ds \{ c_{i} \}\).   
\end{prop}
\textbf{Proof:} 
Starting \(\ds U \in \mySkyB{2+\Omega_{e}}\) and  
 applying the last easy lemma recursively, we have 
 \(\ds U^{(k)} = 0 \) or reach the bottom  
 \( \ds U^{(k)} \in \mySkyL{0} \) where   
 \( U^{(k)}  := 
( \Psi+c_{1}) \circ \cdots \circ  (\Psi+c_{k}) U \)
for some nonzero finite sequence \(\ds (c_{i})\).   
If \( \ds  U^{(k)} = 0\) the discussion finishes.  
If \( \ds  U^{(k)} \ne 0\) we remember a contribution of \(\ds \easyXX{1}{0}\) for the action
of \(\Psi\) and we have \[(\Psi + c_{k+1}) U^{(k)}  = 0 + U^{(k+1)}\;\ 
\text{where}\ 
  U^{(k+1)}\in  \easyXX{1}{0}\;.\]  

Using Lemma \ref{lemma:interest}, we know 
\(\ds \Psi(  U^{(k+1)} ) =   (n+w+1)  U^{(k+1)}\).  
Thus, we get \[ 
( \Psi+c_{1}) \circ \cdots \circ  (\Psi+c_{k}) \circ   (\Psi+c_{k+1})
\circ (\Psi-(n+w+1)) U = 0\;.   \]
\kmqed

\bigskip

\textbf{Proof of Theorem:}
Take \(U\in \myCSW{2}{w,w}\). Combining Propositions
\ref{prop:TL:to:TR}, \ref{prop:TR:to:Rec} and  \ref{prop:Rec:to:Fin},  
we see that 
\[ (\Psi+c_{1}) \circ \cdots \circ (\Psi+c_{m}) U = 0 \] 
for some non-zero numbers \(\ds \{ c_{i} \}\).  
Now express the polynomial  
\( (t +c_{1}) \cdots (t+c_{m}) \) as 
\(  c_{1} \cdots c_{m}  + t g(t) \) for some polynomial \(g(t)\) of one variable \(t\).   

Assume \(U\) is a cycle. Then \(\ds \Psi(U) = (\myHTPmap ) U = (\pdel \Phi ) U\) and 
\begin{align*}
0 &= (\Psi+c_{1}) \circ \cdots \circ (\Psi+c_{m}) U 
=  c_{1} \cdots c_{m} U + \pdel\circ\Phi\circ  g(\pdel\circ\Phi) U \;. 
\end{align*}
This implies \( U\) is exact.  \kmqed

\bigskip

\begin{remark} We can define \(\Phi\) for each m-chain for \(m\geq 2\),
and we may say \(\phi=\Phi\) for 1-chains. 
Thus, it seems to be interesting to study geometry and combinatorics of  
\(\Psi := \pdel \circ \Phi + \Phi \circ \pdel \) for $m$-chains with \(m>2\).
\end{remark}
 
\bibliographystyle{plain}
\bibliography{km_refs}

\def\cprime{$'$} \def\cprime{$'$}
\begin{thebibliography}{1}

\bibitem{Mik:Miz:homogPoisson}
Kentaro Mikami and Tadayoshi Mizutani.
\newblock {Cohomology groups of homogeneous Poisson structures}.
\newblock arXiv:1511.00199v4, May 2017.

\bibitem{Mik:Miz:super2}
Kentaro Mikami and Tadayoshi Mizutani.
\newblock {Euler number and Betti numbers of homology groups of pre Lie
  superalgebra}.
\newblock arXiv:1809.08028v2, December 2018.

\bibitem{Mik:Miz:super}
Kentaro Mikami and Tadayoshi Mizutani.
\newblock {Euler number of homology groups of super Lie algebra}.
\newblock arXiv:1809.08028v1, September 2018.

\end{thebibliography}

\end{document}